\documentclass[12pt]{amsart}

\usepackage{graphicx}

\newcommand{\embfig}[3]{\medskip
\centerline{\includegraphics[width=.#1\textwidth]{#2}}
\centerline{#3}
\medskip
}

\newcommand{\Aaa}{\embfig{80}{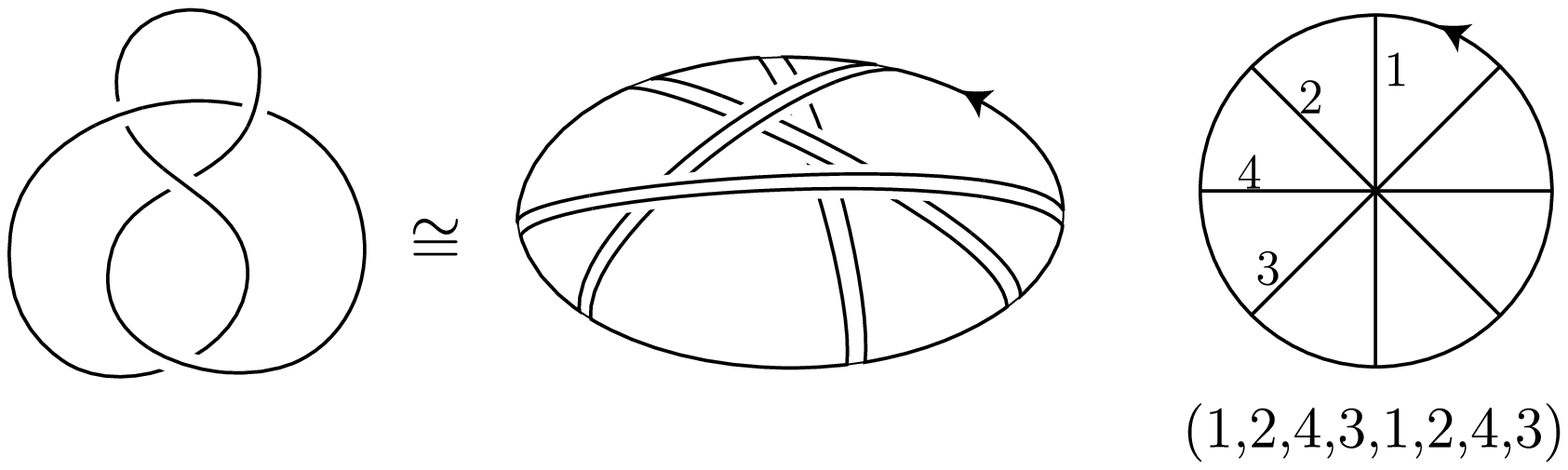}{Figure 1.1:
A flat plumbing basket presentation of 
the knot $4_1$.
}}
\newcommand{\BAa}{\embfig{63}{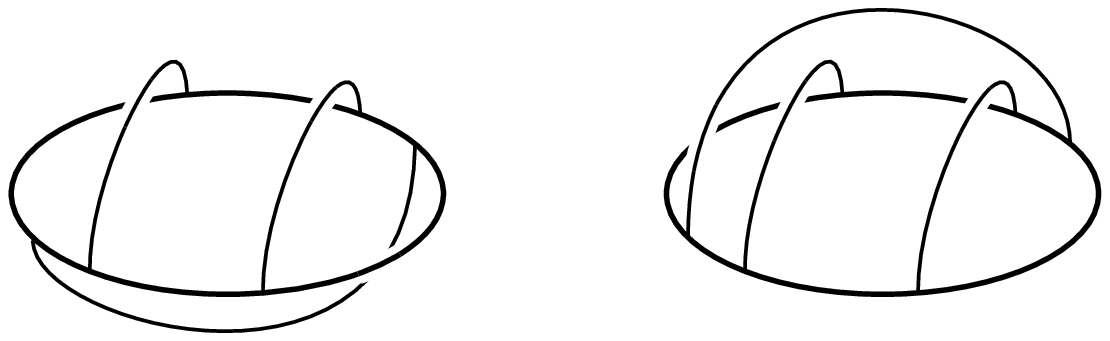}{Figure 2.1:
Flipping a band over.}}
\newcommand{\Baa}{\embfig{90}{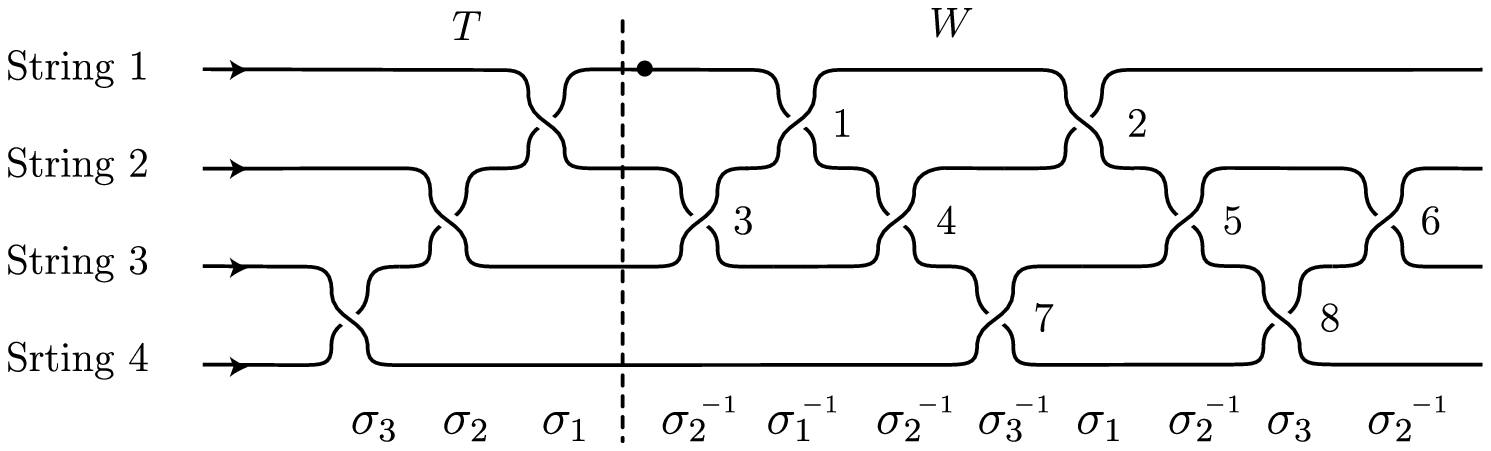}{Figure 2.2:
Labelling the letters in $W$.}}
\newcommand{\Bab}{\embfig{90}{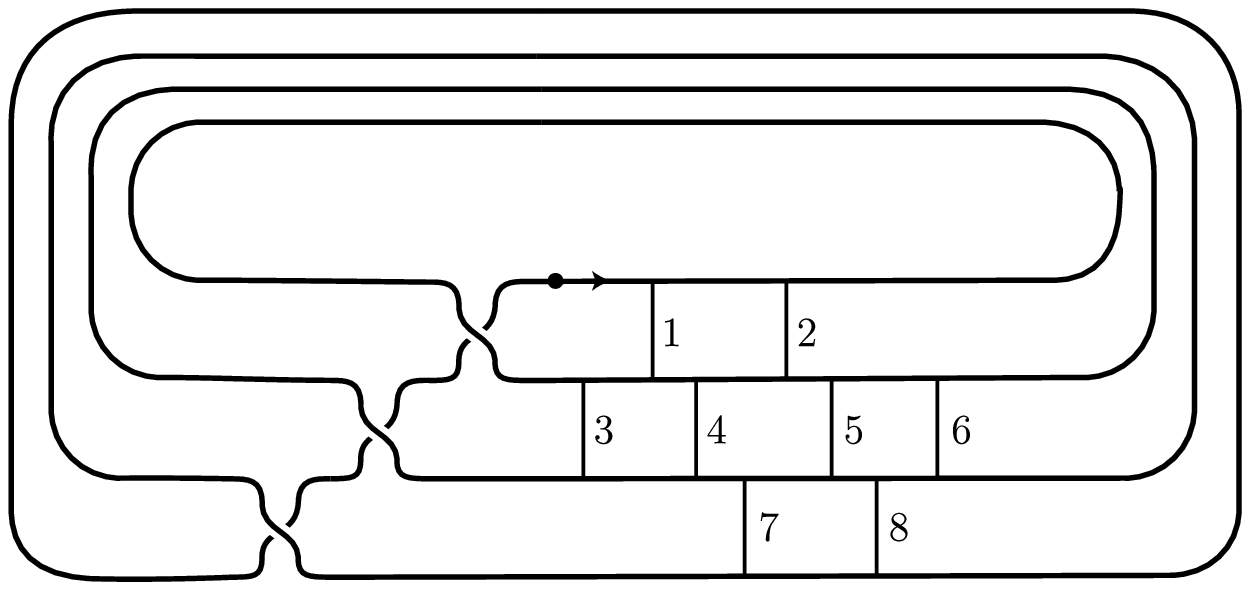}{Figure 2.3:
Reading a flat basket code.}}
\newcommand{\Caa}{\embfig{96}{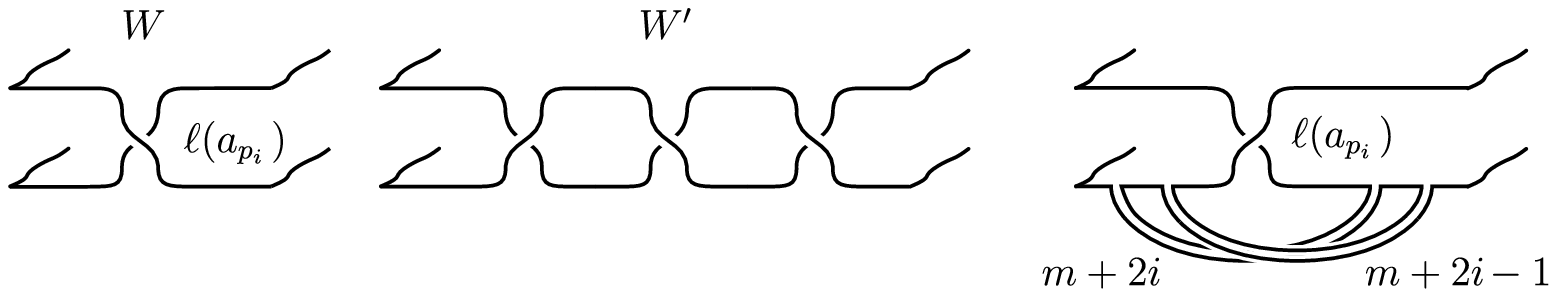}{Figure 3.1:
Deformation of bands in $W'$ and relabeling.
}}
\newcommand{\Cab}{\embfig{95}{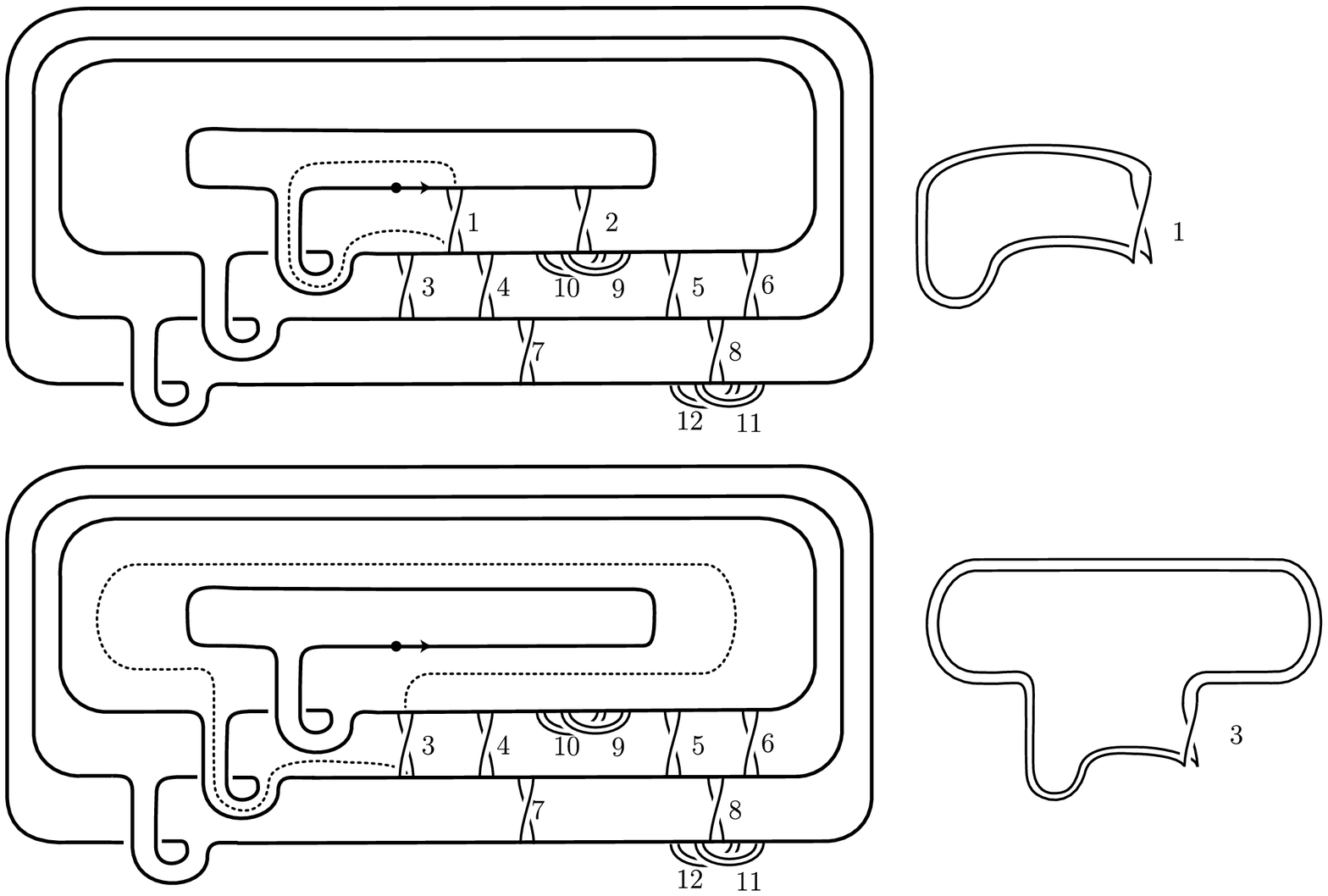}{Figure 3.2:
Removing vertical bands by deplumbing
flat annuli.
}}
\newcommand{\Daa}{\embfig{35}{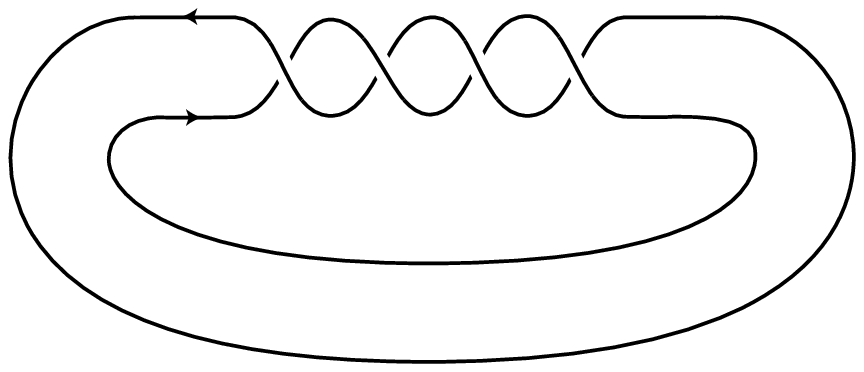}{Figure 4.1:
The oriented $2$-bridge link $S(4,1)$.
}}
\newcommand{\Dbb}{\embfig{60}{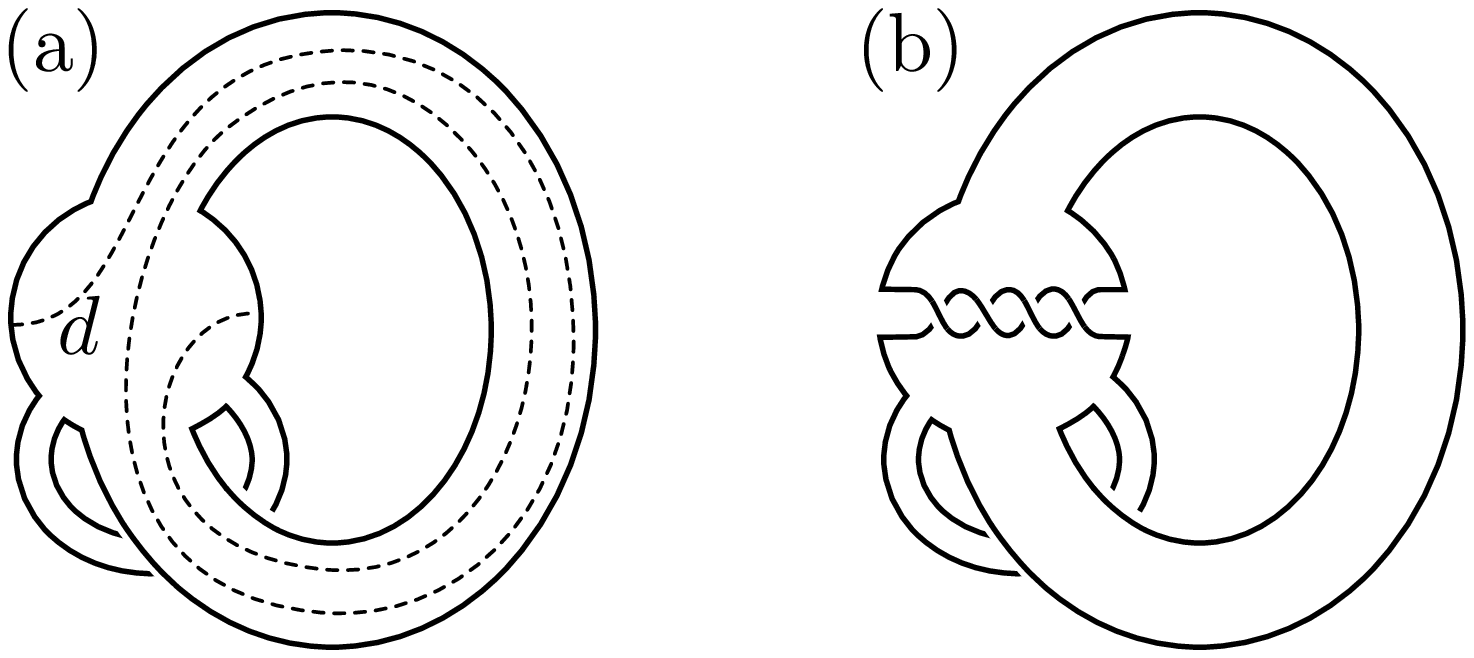}{Figure 4.2:
A flat annuli plumbing for $S(4,1)$ with three
bands
}}
\newcommand{\Dcc}{\embfig{70}{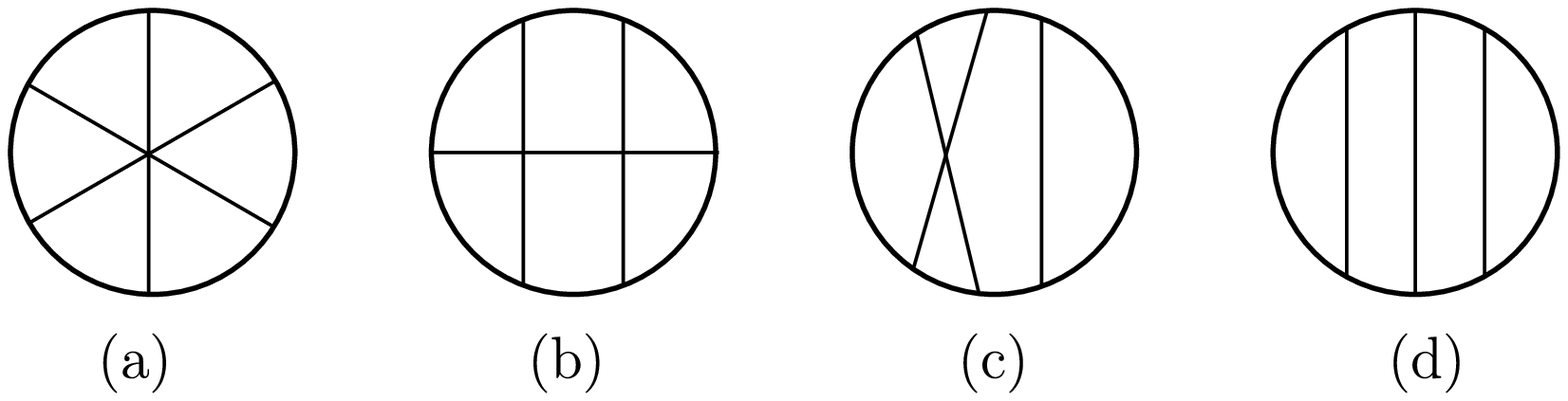}{Figure 4.3:
Flat basket diagrams with three arcs.
}}
\newcommand{\Ddd}{\embfig{32}{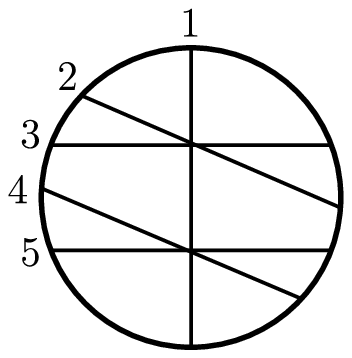}{Figure 4.4:
A flat plumbing basket for $S(4,1)$ with 
five bands.
}}

\newtheorem{thm}{Theorem}[section]

\newtheorem{prop}[thm]{Proposition}

\newtheorem{claim}[thm]{Claim}

 \theoremstyle{definition}
 \newtheorem{dfn}[thm]{Definition}
 \newtheorem{rem}[thm]{Remark}
 \newtheorem{ex}[thm]{Example}

 \title{
Seifert surfaces in open books,
and a new coding algorithm for links}

 \author{Rei Furihata}
 \address{R. Furihata,
 Yosami Junior High School, \newline \indent
5-1 Kamizawatari Ogakie-cho, Kariya city, Aichi 448-0813 Japan}

 \author{Mikami Hirasawa}
\address{M. Hirasawa,
Department of Mathematics,
Gakushuin University, \newline \indent
Mejiro 1-5-1, Toshima-ku Tokyo 171-8588 Japan}
\email{hirasawa@math.gakushuin.ac.jp}

\author{Tsuyoshi Kobayashi}
\address{T. Kobayashi,
Department of Mathematics,\newline \indent
Nara Women's University,
Kitauoya Nishimachi, Nara 630-8506 Japan}
\email{tsuyoshi@cc.nara-wu.ac.jp}

 \date{}
 


\begin{document}
 \maketitle
 
\pagestyle{myheadings}
 \markboth{Seifert surfaces in open books,
 and a new coding algorithm for links}
 {R. Furihata, M. Hirasawa and T. Kobayashi}

\begin{abstract}
We introduce a new standard form of a Seifert surface $F$. In that
standard form, $F$ is obtained by successively plumbing flat annuli
to a disk $D$, where the gluing regions are all in $D$.
We show that any link has a Seifert surface in the standard form,
and thereby present a new way of coding a link.   We present an
algorithm to read the code directly from a braid presentation.
\end{abstract}

\section{Introduction}

The trivial open book decomposition ${\mathcal O}$ of the $3$-sphere $S^3$
is a decomposition of $S^3$ into infinitely many
disk (called pages) sharing their boundaries.
An arc presentation of a link $L$ is a presentation of $L$ as a finite union 
of arcs each of which is properly embedded in a page.
In \cite{Cro}, P. Cromwell introduced the concept of arc presentation
and studied its properties.

It is well known that any link $L$ admits a Seifert surface $F$ (i.e., a compact, orientable, 
embedded surface
whose boundary coincides with $L$).
In this paper, we introduce a concept of positions of Seifert surfaces
called the flat plumbing basket presentation, modifying the definition
of arc presentation.

\begin{dfn}
A Seifert surface is said to be a 
{\it flat plumbing basket}
if it consists of a single page of ${\mathcal O}$ and 
finitely many bands which are
embedded in distinct pages.
 We say that a link $L$ admits a 
 {\it flat plumbing basket presentation} if there exists
 a flat plumbing basket $F$ such that $\partial F$ is equivalent to $L$.
 (See Figure 1.1.)
 \end{dfn}
 
 In this paper, we show that every link admits a flat plumbing basket presentation 
 (Theorem \ref{main} (1)).
  
\Aaa

\begin{rem}
(1)
In \cite{Rud}, L. Rudolph called a Seifert surface $F$ a {\it basket}
if $F$ is obtained from a disk $D$ by  successively plumbing unknotted
(possibly twisted) annuli in such a way that every gluing region is embedded in $D$.
 It is obvious that a flat plumbing basket 
 is obtained from a disk by plumbing flat annuli
in such a way that every gluing region is embedded in the base disk.
Hence a flat plumbing basket is a basket.

(2) In \cite{HaWa}, C. Hayashi and M. Wada showed that any link admits a  Seifert surface
which is obtained from a disk $D$ by successively plumbing flat annuli.
Here we note that the gluing regions in their construction 
are not necessarily contained in $D$.
In Section 4, we give an example showing
 a gap between their \lq flat annuli plumbing' and our
flat plumbing basket in terms of the minimal necessary
number of plumbings to construct a specific link.
\end{rem}

\begin{dfn}
(1) A {\it flat basket diagram} ${\mathcal D}$ is a pair consisting
of a disk $D$ 
(whose boundary has the counter-clockwise orientation)
and a union of 
arcs $a_1, \ldots, a_n$ properly embedded in $D$
with $\partial a_i \cap \partial a_j=\emptyset$\ $(i\neq j)$.
We call the subscript $i$ of $a_i$ 
the {\it label} of the arc.

(2) For a flat basket diagram ${\mathcal D}$, by recording the labels of the arcs as
one travels along $\partial D$, one obtains a word $W$ in $\{1, \ldots, n\}$
such that each letter appears exactly twice.
We call $W$ a {\it flat basket code} for ${\mathcal D}$.
\end{dfn}

\begin{rem}
(See Definition \ref{plumb} for terminology.)
Given a flat plumbing basket $F$, we can eventually obtain a disk by
successively deplumbing flat annuli from $F$ off the positive side of $D$.
Each deplumbing corresponds to removing a band. By labelling the bands of $F$
in the order of these deplumbings, we obtain a flat basket diagram (see Figure 1.1), and hence a 
flat basket code $F$. We also call that flat basket code a  {\it flat basket code}
for the link $L=\partial F$.
Conversely, given a flat basket diagram ${\mathcal D}$, 
we can obtain the flat plumbing basket that ${\mathcal D}$ presents by attaching
flat bands $b_n, b_{n-1}, \ldots, b_1$ on the positive side
along $a_n, a_{n-1}, \ldots, a_1$ in this order.
Here we note that $b_i$ is closer to $D$ than $b_{i+1}$. 
\end{rem}

\begin{ex}
The flat basket code $(1,2,4,3,1,2,4,3)$ represents
the figure-eight knot $4_1$ (see Figure 1.1), 
and the code
$(1,2,3,4,1,2,3,4)$ represents the trefoil $3_1$.
We see that any flat plumbing basket for
$4_1$ or $3_1$ has at least four bands, since 
these two knots have genus 1, and any flat plumbing
basket is compressible and hence not of minimal genus. 
\end{ex}

In Theorem \ref{main} (2), 
we give an algorithm to obtain a flat basket code for any link $L$,  
using a braid presentation of $L$.

\begin{rem}
Let ${\mathcal C}_L$ be the set
of flat basket codes for a link $L$.
Each ${\mathcal C}_L$ is non-empty, and has
a unique element $c_L$ which has the shortest length
and lexicographically minimal.
Therefore, we can regard the set of all link types 
as a totally ordered set, by defining $L_i < L_j$
if and only if $c_{L_i}<c_{L_j}$.
\end{rem}

\noindent
Acknowledgements.
The authors would like to express their thanks to 
Prof. Kazuhio Ichihara for helpful comments.

\section{Preliminaries and the statement of the
results}
First we recall the notion of {\it plumbing} and
{\it deplumbing} of Seifert surfaces.

\begin{dfn}\label{plumb}
Let $R$ be a Seifert surface in $S^3$. We say that $R$ is a plumbing of surfaces
$R_1$ and $R_2 (\subset S^3)$ if the following conditions are satisfied:

(1) $R=R_1\cup R_2$ and $R_1 \cap R_2=D$, 
where $D$ is a square with four successive edges
$a_1, b_1, a_2, b_2$ such that $a_i$ (resp. $b_i$) is contained in 
$\partial R_2$
(resp. $\partial R_1$) and properly embedded in $R_1$ (resp. $R_2$).\\
(2) There exist $3$-balls $B_1, B_2 \subset S^3$ satisfying the following:\\
\indent
(2.1) $B_1\cup B_2=S^3$ and $B_1\cap B_2=\partial B_1=\partial B_2$.\\
\indent
(2.2) $B_i \supset R_i, \ i=1, 2$.\\
\indent
(2.3) $\partial B_1\cap R_1 = \partial B_2\cap R_2=D$.

Here the square $D$ is called the {\it gluing region} of the plumbing.
We may say that
$R$ is obtained from $R_1$ by plumbing $R_2$ along $D$,
or that $R_1$ is obtained from $R$ by deplumbing $R_2$.
\end{dfn}

\begin{dfn}
Suppose a Seifert surface
$R$ is a plumbing of $R_1$ and $R_2$ along the gluing region $D$.

(1)
Let $T$ be a subsurface of $R_1$. We say that $R_2$ is {\it rooted in}
$T$, if $D$ is contained in $T$.

(2)
Since $R_1$ is oriented, it admits two sides, say the $+$side and $-$side.
If the $3$-ball $B_2$ lies in the $\varepsilon$-side 
($\varepsilon = \pm $), we say that $R$ is obtained from $R_1$ by plumbing $R_2$ 
along $D$ on the $\varepsilon$-side.
Conversely, we say that $R_1$ is obtained from $R$ by deplumbing
$R_2$ off the $\varepsilon$-side.

(3)
Let $a$ be an arc properly embedded in $R_1$,
Let $N$, be a regular neighborhood of $a$ in $R_1$.
Push $N-\partial R_1$ into the 
$\varepsilon$-side slightly to make a band $N'$ attached to $R_1$.
Let $R=R_1 \cup N'$. Then $R$ is a plumbing of $R_1$ and a flat annulus
(in fact, $N\cup N'$) on the $\varepsilon$-side.
We say that $R$ is obtained from $R_1$ by attaching a flat band  along
$a$ on the $\varepsilon$-side.

\end{dfn}

\begin{rem}
(1) A Seifert surface $F$ is isotopic to a flat plumbing basket
if and only if
$F$ is obtained from a disk $D$ by successively plumbing flat
annuli each of which is rooted in $D$.\\
(2) In (1), some of the annuli may be plumbed on the $-$side.
However, by using isotopy as in Figure 2.1, we can regard all
the annuli as plumbed on the $+$side.
\end{rem}

\BAa

Note that any link can be expressed as the closure of some
$n$-braid \cite{vo},
in particular, with a braid word of the form
$\sigma_{n-1}\sigma_{n-2}\cdots\sigma_{1} W$.
(In fact, if $L$ is the closure of a braid $B$, we may let 
$W=(\sigma_{n-1}\sigma_{n-2}\cdots\sigma_{1})^{-1} B$.)
We adopt the following
 convention on braids.
 
 \noindent
{\bf Convention on braids.} See Figure 2.2.
In drawing braids, we put strings horizontally and
call the top string the first string and so forth.
Strings are oriented from left to right,
and we read the braid word from the left end.
The generator $\sigma_i$ corresponds to a crossing
of the $i^{\rm th}$ and $(i+1)^{\rm st}$ 
strings, where the $i^{\rm th}$ string goes down crossing
over the $(i+1)^{\rm st}$ which goes up.

The main result of this paper is as follows:

\begin{thm}\label{main}
Let $L$ be an oriented link.
Express $L$ as a closed $n$-braid with a braid word
$\sigma_{n-1}\sigma_{n-2}\cdots\sigma_{1} W$.
Suppose the length of $W$ is $m$ and that $W$ has $s$ positive letters.
Then we have the following:

(1) There exists a flat plumbing basket $F$ with $m+2s$ bands
such that $\partial F$ is isotopic to $L$.

(2)
A flat basket code ${\mathcal C}$ 
(of length $2(m+2s)$) coming from $F$ is
obtained by the following algorithm.

\noindent
{\bf Algorithm.} Let $a_i$ be the $i^{\rm th}$ letter in $W$.

Step 1:
We give distinct labels in $\{ 1, \dots , m\}$ to $a_i$'s.
The label $\ell(a_i)$
is defined to be $k$ if $a_i$ is in the $k^{\rm th}$ position
with respect to the lexicographic order by the double index of $a_i$,
where the double index of $a_i$ is $(q,i)$ if $a_i=\sigma_q^{\pm 1}$.

For each $i =1, \ldots, n$, let
$M_i$ denote the word obtained from the word $a_1 a_2\cdots a_m$ as follows:\\
$M_1$ is obtained 
by deleting all $a_j$'s which are not  $\sigma_1^{\pm1}$.\\
$M_n$ is obtained
by deleting all $a_j$'s which are not  $\sigma_{n-1}^{\pm1}$.\\
$M_i\  (1<i<n)$ is obtained 
by deleting all $a_j$'s which are not  $\sigma_{i-1}^{\pm1}$ or 
$\sigma_{i}^{\pm1}$.

Let ${\mathcal C}_1$ be the word obtained from the composed word
$M_1 M_2 \cdots M_n$ by replacing each $a_j$ by $\ell(a_j)$.
Note that ${\mathcal C}_1$ is a word in 
$\{1, 2, \ldots, m\}$ such that each letter appears exactly twice.

Step 2:
Let $\{a_{p_1}, \ldots, a_{p_s}\}\  (p_1<p_1<\cdots<p_s)$ be the subset of 
$\{a_1, \ldots, a_m\}$ such that each element corresponds to a positive letter.
Then the desired code ${\mathcal C}$ is obtained as follows:\\
In ${\it \mathcal C}_1$, replace the second appearance
of each $\ell(a_{p_i})$, $i=1, \ldots, s$, by a sequence 
$m+2i, m+2i-1, \ell(a_{p_i}),
m+2i, m+2i-1$.
\end{thm}

 \Baa

  \begin{ex}\label{ex:1}
  See Figure 2.2. Here, $W=
 \sigma_2^{-1}\sigma_1^{-1}\sigma_2^{-1}
 \sigma_3^{-1}\sigma_1^{} \ \sigma_2^{-1}\sigma_{3}^{ }\ 
  \sigma_2^{-1}$ and $m=8, s=2$.
The double indices for $a_1, \ldots, a_8$ are
$(2,1)$, $(1, 2)$, $(2,3)$, $(3,4)$, $(1,5)$, $(2,6)$,
$(3,7)$, $(2,8)$.
Since
$(1,2)<(1,5)<(2,1)<(2,3)<(2,6)<(2,8)<(3,4)<(3,7)$,
we have:\\
 $(\ell(a_1),\ell(a_2),$
 $\ell(a_3),\ell(a_4),
 \ell(a_5),\ell(a_6),$
 $\ell(a_7),\ell(a_8))=$
 $(3,1,4,7,2,5,8,6)$.
In Figure 2.2,  numbers put beside the crossings present the labels.
On the other hand, we have 
$M_1=a_2a_5,\
M_2=a_1a_2a_3a_5a_6a_8$,
$M_3=a_1a_3a_4a_6a_7a_8,\
M_4=a_4a_7$.
Hence we have:
${\mathcal C}_1=
\ell(a_2)\ell(a_5)\cdot\\
\ell(a_1)\ell(a_2)\ell(a_3)\ell(a_5)\ell(a_6)\ell(a_8)\cdot
\ell(a_1)\ell(a_3)\ell(a_4)\ell(a_6)\ell(a_7)\ell(a_8)\cdot
\ell(a_4)\ell(a_7)$
$=$
$(1\ 2\ 3\ 1\ 4\ 2\ 5\ 
6\ 3\ 4\ 7\ 5\ 8\ 6\ 7\ 8)$.
Since $a_5$ (with label $2$) and $a_7$ (with label
$8$) correspond to the positive letters,
we replace the second appearances of $2$ and $8$
in ${\mathcal C}_1$ respectively by
$(10 \ 9\ 2\ 10\ 9)$ and $(12\ 11\ 8\ 12\ 11)$.
Thus, we obtain
a flat basket code 
$${\mathcal C}=
(1\ 2\ 3\ 1\ 
4\ 10 \ 9\ 2\ 
10\ 9\ 5\ 
6\ 3\ 
4\ 7\ 5\ 8\ 6\ 
7\ 12\ 11\ 8\ 
12\ 11).$$
Note that we can read off ${\mathcal C}_1$, and hence ${\mathcal C}$
from a geometric braid.
See Figure 2.3, which is a schematic picture
where the  crossings corresponding to the letters of
$W$ are represented by arcs. The arcs are labelled
from \lq top-left\rq \ to \lq bottom-right\rq.
In this picture, the closed braid $B$ obtained
by removing all the crossings corresponding to the
letters in $W$ is the unknot.
We obtain the code ${\mathcal C}_1$ by travelling
along $B$ with the starting point depicted by the dot
and recording the labels of the arcs as we pass their endpoints.
 \end{ex}

 \Bab

\section{Proof of Theorem 2.4.}

Suppose that
$L$ is the closure of an $n$-braid with a braid word $TW$, 
where $T$ is $\sigma_{n-1}\cdots \sigma_2\sigma_1$.
Let $W'$ be the word obtained from $W$ by replacing each positive
letter $\sigma_i$ by $\sigma_i^{-1}\sigma_i \sigma_i$.
Note that $L$ is also presented by $TW'$.
Let $F$ be the Seifert surface for $L$ obtained from the closure of 
the geometric braid presented by $TW'$ by applying Seifert algorithm, i.e.,
$F$ is obtained from $n$ parallel horizontal disks 
$D_1, \ldots, D_n$
by attaching half twisted bands each corresponding to a letter in $TW'$,
where $D_i$ corresponds to the $i^{\rm th}$ string and, as we may suppose,
the heights of $D_1, \ldots, D_n$ are in this order with $D_1$ the highest.
We assume the face-side of $D_1$ coincides with
the $+$side of $F$.

Let $D$ be the disk in $F$ consisting of
$D_1, \ldots, D_n$  and 
the $n-1$ bands corresponding to the letters in
$T=\sigma_{n-1}\cdots \sigma_2\sigma_1$. From now on, 
we regard $F$ as consisting of
the disk $D$ and bands attached to $D$.
Namely, $T$ corresponds to the disk $D$,
and $W'$ corresponds to the bands attached to $D$.

Note that each band from $W$ corresponds to some $a_j$.
Then label the band with $\ell(a_j)$.
Recall $\{a_{p_1},\ldots, a_{p_s}\}$ in Step 2.
As we obtain $F$ from the Seifert surface obtained from $TW$, 
we locally attach two bands near each band
labelled $\ell(a_{p_i})$. See Figure 3.1.
Deform two of
the bands as specified in Figure 3.1 and (re)label
the three bands as $m+2i-1, m+2i, \ell(a_{p_i})$. 

\Caa

Now there are $m+2s$ bands attached to $D$.
We call the first $m$ bands
(labelled from $1$ to $m$) the {\it vertical bands},
since each of them connects two distinct sub-disks
$D_t$ and $D_{t+1}$ of $D$.
On the other hand, we call the last $2s$ bands (labelled from $m+1$ to $m+2s$) the
{\it horizontal bands} since each of them connects a single sub-disk of $D$.
Note that all the vertical bands are twisted in the same way. 
Actually, the system of vertical bands
corresponds to the word obtained from $W$ by
replacing all $\sigma_{*}$ by $\sigma_{*}^{-1}$.
We denote this Seifert surface for $L$ by $F'$ and
regard the disk $D$ as a subsurface in $F'$.

We prove the following:

\begin{claim}
The surface $F'$ thus obtained
is isotopic to a flat plumbing basket for $L$,
where all bands are rooted in $D$.
\end{claim}

We show that we can deplumb $m+2s$ flat annuli, say,
$A_1, A_2,\ldots, A_{m+2s}$ each rooted in $D$ successively in this order,
in such a way that deplumbing $A_i$ corresponds to removing the band labelled $i$.

Let $h_1, h_2, \ldots, h_{n-2}, h_{n-1}=m$ be integers such that
(i) there exist $h_1$ vertical bands connecting $D_1$ and $D_2$,
(ii) there exist $h_i - h_{i-1}\  (2\le i\le n-1)$ vertical bands connecting
$D_i$ and $D_{i+1}$. 

See Figure 3.2 for the same braid used in Example \ref{ex:1},
where the bands connecting the subdisks of the base disk $D$
are deformed so that only one side of $D$ is visible.
Let $F_1$ be the surface obtained from $F'$ by removing the bands
labelled $1, \ldots, h_1$. 
It is directly observed from Figure 3.2 that 
$F_1$ is obtained from $F'$ by deplumbing a sequence of flat annuli,
say,  $A_1, \ldots, A_{h_1}$ (where $A_j$ corresponds to the band $j$)
each rooted in $D$. Here, each of the gluing regions 
(represented by dotted lines) goes once though
the band in $D$ connecting $D_1$ and $D_2$.
Then repeat the same procedure for the bands 
labelled $i\  (h_1<i \le h_2)$
to obtain a new surface $F_2$ and so forth.

\Cab

Finally we obtain a surface $F_{n-1}$ which is a union of $D$ and
the horizontal bands.
We can remove the $2s$ horizontal bands by deplumbing flat annuli rooted in $D$,
say, $A_{m+1},\ldots, A_{m+2s}$ in this order,
where $A_j$ corresponds to the band $j$.
Finally, we have constructed a Seifert surface for $L$ which is isotopic
to a flat plumbing basket with $m+2s$ bands.
The proof of Theorem \ref{main} (1) is now completed.

For the proof of Theorem \ref{main} (2), it is enough to confirm that a flat basket
code of $F'$ obtained by using $A_1, \ldots, A_{m+2s}$ coincides with the code
that the algorithm yields. This can be done by travelling along $\partial D$
starting from the dot indicated in Figure 2.3.

\section{Examples}
In \cite{HaWa}, Hayashi and Wada showed that any oriented link
admits a Seifert surface which is obtained from a disk by plumbing
flat annuli. A flat plumbing basket can be regarded as such an example,
but not vice versa.
In this section, we exhibit a Seifert surface which is obtained from
a disk by successively plumbing flat annuli, but which is not isotopic to 
any flat plumbing basket.

\begin{prop} \label{difference}
Let $L$ be the $2$-bridge link $S(4,1)$,
 oriented so that $L$ spans an unknotted
annulus with two full-twists (Figure 4.1). 
Then we have the following.\\
(1) $L$ admits a Seifert surface which is obtained from a disk by plumbing three
flat annuli.\\
(2) Any flat plumbing basket for $L$ has at least $5$ bands.
\end{prop}

\Daa

{\it Proof.}
Figure 4.2 shows Proposition \ref{difference} (1):
Consider a dotted line $d$ in a surface obtained by
plumbing two flat annuli to a disk, as in Figure 4.2 (a).
Then by plumbing another flat annulus along $d$
on the face-side, 
we obtain a Seifert surface
isotopic to that in Figure 4.2 (b), whose boundary is equivalent to $L$.
Note that the gluing region for the third plumbing goes through
a previously plumbed flat annulus.

\Dbb

To prove (2), we examine flat plumbing baskets which have
$\le 4$ bands. Suppose a flat plumbing basket $F$ has
$n$ bands.
If $n\le2$, then it is easy to observe that $\partial F$ is either the trivial knot or
a trivial link with two or three components.
If $n =3$, then there are four possible patterns of ambient arcs of flat basket 
diagrams depicted in Figure 4.3.

\Dcc

In Figure 4.3 (a), $\partial F$ is a (positive or negative) Hopf link.
In other cases, $\partial F$ is a trivial link with two or four components.
If $n=4$, then $\partial F \neq L$ since $\partial F$ has an odd number of
connected components. (Note that
removing a band changes the number of connected components
by $\pm 1$).
Hence Proposition \ref{difference} is proved.
\qed

Note that 
$S(4,1)$ admits a flat plumbing basket 
presentation with
five bands as in
Figure 4.4, with a flat basket code $(1,2,3,4,5,1,4,5,2,3)$.

\Ddd

\end{document}